\documentclass{article}
\usepackage{amsmath,epsf}
\usepackage[latin1]{inputenc}
\usepackage{amsfonts}

\title{On certain spaces of \\lattice diagram determinants}
\author{Jean-Christophe Aval}
\date{Laboratoire A2X, Universit\'e Bordeaux 1\\ 351 cours de la Lib\'eration, F-33405 Talence cedex\\ e-mail : {\tt aval@math.u-bordeaux.fr}}

\numberwithin{equation}{section}
\newtheorem{thm}{Theorem}[section]
\newtheorem{prop}[thm]{Proposition}
\newtheorem{conj}[thm]{Conjecture}
\newtheorem{lem}[thm]{Lemma}

\newtheorem{definition}[thm]{Definition}

\newtheorem{rem}[thm]{Remark}

\font \sc=cmr9

\def\N{{\mathbb N}}

\def\Q{{\mathbb Q}}
\def\S{{\mathcal S}}
\def\M{{\mathcal M}}
\def\T{{\mathcal T}}
\def\I{{\mathcal I}}
\def\F{{\mathcal F}}
\def\DD{\Delta}

\def\Young#1{\vbox{\smallskip\offinterlineskip
    \halign{&\vbox{##}\kern-\Thickness\cr #1}}}

\newdimen\Squaresize \Squaresize=20pt
\newdimen\Thickness \Thickness=.1pt
\newdimen\Correction \Correction=7pt

\def\Vide#1{\hbox{
       \vbox to \Squaresize{\vss
          \hbox to \Squaresize{\hss#1 \hss}\vss}
    \hskip-\Correction}
   \kern-\Thickness}

\def\Carre#1{\hbox{\vrule width \Thickness
   \vbox to \Squaresize{\hrule height \Thickness\vss
      \hbox to \Squaresize{\hss$\scriptstyle#1$\hss}
   \vss\hrule height\Thickness}
   \unskip\vrule width \Thickness}
   \kern-\Thickness}

\begin{document}

\maketitle

\begin{abstract}
The aim of this work is to study some lattice diagram determinants $\Delta_L(X,Y)$ as defined in \cite{lattice} and to extend results of \cite{untrou}. We recall that $M_L$ denotes the space of all partial derivatives of $\Delta_L$. In this paper, we want to study the space $M^k_{i,j}(X,Y)$ which is defined as the sum of $M_L$ spaces where the lattice diagrams $L$ are obtained by removing $k$ cells from a given partition, these cells being in the ``shadow'' of a given cell $(i,j)$ in a fixed Ferrers diagram. We obtain an upper bound for the dimension of the resulting space $M^k_{i,j}(X,Y)$, that we conjecture to be optimal. This dimension is a multiple of $n!$ and thus we obtain a generalization of the $n!$ conjecture. Moreover, these upper bounds associated to nice properties of some special symmetric differential operators (the ``shift'' operators) allow us to construct explicit bases in the case of one set of variables, i.e. for the subspace $M^k_{i,j}(X)$ consisting of elements of $0$ $Y$-degree. 
\end{abstract}

\section{Introduction}

\begin{definition}\label{mudef}
A {\sl lattice diagram} is a finite subset of $\N\times\N$.
For $\mu_1\geq \mu_2\geq \cdots \geq \mu_k>0$, we say that $\mu=(\mu_1,
\mu_2, \dots,
\mu_k)$ is a {\sl partition} of $n$ if $|\mu|=\mu_1+\cdots +\mu_k$ equals $n$. We associate to
a partition $\mu$ its Ferrers diagram $\{(i,j)\,:\,0\le i\le k-1,\,
0\le j\le \mu_{i+1}-1\}$ and we shall use the symbol
$\mu$ for both the partition and its Ferrers diagram.
\end{definition}
Most definitions and conventions we use are similar to \cite{lattice}.
For example, given the partition $(4,2,1)$, its partition diagram is
  $$\Young{\Carre{2,0}\cr
           \Carre{1,0}&\Carre{1,1}\cr
           \Carre{0,0}&\Carre{0,1}&\Carre{0,2}&\Carre{0,3}\cr
           }\quad.$$
It consists of the lattice cells
$\{(0,0),(1,0),(2,0),(0,1),(1,1),(0,2),(0,3)\}$.

Let now $X=X_m=\{x_1,x_2,\dots,x_m\}$ and $Y=Y_m=\{y_1,y_2,\dots,y_m\}$ be two sets of $m$ variables and $\Q[X]\!=\!Q[x_1,x_2,\dots,x_m]$ and $\Q[X,Y]\!=\!Q[x_1,x_2,\dots,x_m,\\ y_1,y_2,\dots,y_m]$ denote respectively the rings of polynomials in $m$ and $2m$ variables with rational coefficients. Since we have to deal with polynomials in $\Q[X]$ or $\Q[X,Y]$, we shall denote by $Z$ a subalphabet of $(X,Y)$ and by $\Q[Z]$ the corresponding ring of polynomials.

\begin{definition}
Given a lattice diagram $L=\{(p_1,q_1), (p_2,q_2),\dots , (p_n,q_n)\}$ with $n$ cells
we define the {\sl lattice determinant}
\begin{equation}
  \Delta_L(X,Y)= \det \big( x_i^{p_j}y_i^{q_j}\big)_{1\le i,j\le n}.
\end{equation}
\end{definition}
The polynomial
$\Delta_L(X,Y)\in\Q[X_n,Y_n]=\Q[X,Y]$ (with $m=n$, the number of cells in the diagram $L$) is different from zero only if the diagram $L$ consists of $n$ distinct cells in the positive quadrant. In this case $\DD_L$ is bihomogeneous of degree $|p|=p_1+\cdots +p_n$ in $X$
and of degree $|q|=q_1+\cdots +q_n$  in  $Y$. To insure that this definition
associates a unique determinant to $L$ we require
that the list of lattice cells is given with respect to the lexicographic order with priority to the second entry that is to say:
\begin{equation}\label{lex}
(p_1,q_1)<(p_2,q_2) \quad\iff\quad q_1<q_2\quad\hbox{or}\quad [q_1=q_2\hbox{ and }p_1<p_2].
\end{equation}

\begin{definition}
For a polynomial $P(Z)\in\Q[Z]$, the vector space spanned by all the partial
derivatives of $P$ of
all orders is denoted
${\mathcal L}_\partial[P]$, i.e.:
\begin{equation}
{\mathcal L}_\partial[P]=\Q[\partial Z],
\end{equation} 
where for a polynomial $Q$ in $\Q[Z]$, $Q(\partial)=Q(\partial Z)$ denotes the differential operator obtained by substituting $x_i$ and $y_i$ respectively by $\partial x_i$ and $\partial y_i$ in the expression of $Q$.
Next we define
\begin{equation}  
M_L = {\mathcal L}_\partial[\Delta_L(X,Y)]
\end{equation}
the vector space associated to the lattice diagram $L$.
\end{definition}
A permutation $\sigma\in \S_n$ acts diagonally on
a polynomial $P(X,Y)\in\Q[X_n,Y_n]$ as follows:
$\sigma  P(X,Y)\,=\, P(x_{\sigma_1},x_{\sigma_2},\dots
,x_{\sigma_n},y_{\sigma_1},y_{\sigma_2},\dots ,y_{\sigma_n})$.
Under this action, $\Delta_L(X,Y)$ is clearly an alternant. It follows that for
any lattice diagram $L$ with $n$ cells, the vector space
$M_L$
is an $\S_n$-module. Since $\Delta_L(X,Y)$ is  bihomogeneous, this module
affords a natural
bigrading. 

The most general problem considered in \cite{lattice} and \cite{coltrouee} concerns the space $M_L$. The main question is to decide whether this space is $\S_n$-isomorphic to a sum of left regular representations or not. In \cite{coltrouee}, the case where all the lattice cells of $L$ lies on a single axis is solved. In the particular case where $L$ corresponds to a partition $\mu$ the question leads to the ``$n!$ conjecture'' which asserts that the space $M_\mu$ is a single copy of the left regular representation. Many efforts to prove this conjecture were only sufficient to obtain it in some special cases (see \cite{allen}, \cite{aval}, \cite{gh}, \cite{orbit} for example).

The next class of lattice diagrams that is of interest is obtained by removing a single cell from a partition diagram. Its interest comes in part from the fact that it gives a possible recursive approach for the $n!$ conjecture, with the statement of a conjectural ``four term recurrence.''
If $\mu$ is a partition of $n+1$, we  denote by $\mu/ij$ the lattice
diagram obtained by
removing the cell $(i,j)$ from the Ferrers diagram of $\mu$. We  refer to the cell
$(i,j)$ as
the {\sl hole} of  $\mu/ij$. It is conjectured in \cite{lattice} that the number of copies of the left regular
representations in $M_{\mu/ij}$ is equal to the cardinality (which we denote by $s_{\mu}(i,j)$ or by $s$ if there is no ambiguity) of the
$(i,j)$-{\sl shadow}, where the shadow of a cell $(i,j)$, as shown on the figure below is:
$S_{\mu}((i,j))=\{(i',j')\in\mu : i'\ge i,\,j'\ge j\}$. 

\begin{figure}[h]
$$
\epsfbox{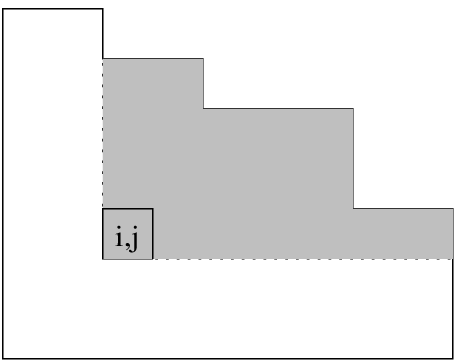}
$$
\end{figure}

A study of the subspace $M_{\mu/ij}(X)$ of $M_{\mu/ij}$ consisting of elements of 0 $Y$-degree can be found in \cite{untrou}, in which the corresponding ``four term recursion'' is proven by using the construction of explicit bases.

The aim of this article is to propose a generalization for the $n!$ conjecture. The space that we consider is defined as follows. Let $\mu$ be a partition of $n+k$. This partition is fixed and does not appear in the following notations. 

\begin{definition}
Let $M_{i,j}^k$ denote the following sum of vector spaces
\begin{equation}\label{sumdef}
M_{i,j}^k=M_{i,j}^k(X,Y)=\sum_{(a_1,b_1),\dots,(a_k,b_k)}M_{\mu/\{(a_1,b_1),\dots,(a_k,b_k)\}},
\end{equation}
where the sum is over all the $k$-tuples of cells in the shadow of $(i,j)$. 
\end{definition}
We first observe that because of the ``shift'' operators (see \cite{lattice}, Proposition I.3 or Section 2 in this paper) we have $M_{\mu/ij}=M_{i,j}^1$ (equation \ref{onehole}). 
Hence this space $M_{i,j}^k$ is a possible generalization of $M_{\mu/ij}$ if we want to make $k$ holes in the Ferrers diagram.
The object of this paper is to show the interest of the space $M_{i,j}^k$ and to give support to the Conjecture \ref{conje} that $\dim M^k_{i,j} = {s \choose k} n!$.

The organization of the article is the following. In the second section we introduce some ``shift'' operators which are useful to move the holes and the cells in the diagrams. The third section is devoted to the proof of an upper bound (${s \choose k} n!$) for the dimension of $M_{i,j}^k$, that is conjectured to be optimal. In the fourth section we study $M_{i,j}^k(X)$, the subspace of $M_{i,j}^k(X,Y)$ consisting of elements of 0 $Y$-degree, for which we obtain explicit bases. 

\section{The ``shift'' operators}

In this paragraph, we want to describe the action of some special symmetric differential operators on the determinants $\Delta_{L}$. We recall the following definitions as stated in \cite{macdonald}: 
\begin{definition}
For each integer $r\ge 1$, the $r$-th power sum $P_r(X)$ (we do not use the classical notation $p_r$ to avoid a possible confusion with the biexponents) is defined by
\begin{equation}
P_r(X)=\sum x_i^r.
\end{equation}
For each integer $r\ge 0$, the $r$-th elementary symmetric function $e_r(X)$ is the sum of all products of $r$ distinct variables $x_i$, so that $e_0=1$ and for $r\ge 1$:
\begin{equation}
e_r(X)=\sum_{i_1<\cdots<i_r}x_{i_1}\cdots x_{i_r}.
\end{equation}
For each integer $r\ge 0$, the $r$-th complete symmetric function $h_r(X)$ is the sum of all monomials of total degree $r$ in the variables $x_i$, so that:
\begin{equation}
h_r(X)=\sum_{i_1\le\cdots\le i_r}x_{i_1}\cdots x_{i_r}.
\end{equation}
\end{definition}
For the sake of simplicity, we only state the following propositions for $X$-shifts. Of course similar results also hold for $Y$-shifts. The only difference concerns the signs. The choice of the lexicographic order \ref{lex} is made to simplify the results and the proofs for $X$-shifts.

\begin{prop}\label{proP}
Let $L$ be a lattice diagram. Then for any integer $k\ge 1$ we have
\begin{equation}\label{pkop}
P_k(\partial X)\Delta_L(X,Y)=\sum_{i=1}^n\pm\epsilon(L,P_k(i;L))\Delta_{P_k(i;L)}(X,Y),
\end{equation}
where ${P_k(i;L)}$ is the diagram obtained by replacing the $i$-th biexponent $(p_i,q_i)$ by $(p_i-k,q_i)$ and the coefficient $\epsilon(L,P_k(i;L))$ is a positive integer. The sign in \ref{pkop} is the sign of the permutation that reorders the obtained biexponents with respect to the lexicographic order \ref{lex}.
\end{prop}

\noindent{\it Proof.}
This is a particular case of Proposition I.1 in \cite{lattice}, but we shall give here a simple proof because some ingredients will be useful later.

If the diagram $L$ consists of the cells $L=\{(p_1,q_1),\dots,(p_n,q_n)\}$, we can develop the determinant $\DD_L$  with respect to the $j$-th column and write:
\begin{equation}\label{dv}
\DD_L(X,Y)=\sum_{i=1}^n x_j^{p_i}y_j^{q_i}.A_{i,j}
\end{equation}
where $A_{i,j}$ denotes the cofactor $(i,j)$. Let us remark that this cofactor is a polynomial where the variable $x_j$ does not appear. Thus when we derive \ref{dv}, we obtain:
\begin{equation}\label{paeq}
\partial x_j^k \DD_L(X,Y)=\sum_{i=1}^n c_i^k x_j^{p_i-k}y_j^{q_i}.A_{i,j}
\end{equation}
where $c_i^k=p_i(p_i-1)\cdots(p_i-k+1)$.
Next we sum \ref{paeq} over $j$ to get:
\begin{equation}\label{som}
P_k(\partial X) \DD_L(X,Y)=\sum_{i=1}^n c_i^k \sum_{j=1}^n x_j^{p_i-k}y_j^{q_i}.A_{i,j}.
\end{equation}
Thus we obtain \ref{pkop} by recognizing in \ref{som} the development (up to sign) of $\Delta_{P_k(i;L)}$. As a biproduct we observe that $\epsilon(L,P_k(i;L))=c_i^k=p_i(p_i-1)\cdots(p_i-k+1)$ and that this coefficient does not depend on the operator $P_k$.

\begin{rem}\rm
The diagram ${P_k(i;L)}$ is the diagram obtained by pushing down the $i$-th cell of $L$: its biexponent $(p_i,q_i)$ is replaced by $(p_i-k,q_i)$ which corresponds to $k$ steps down. The other biexponents are unchanged. This duality between the substractions on the set of biexponents and the movements of cells in the diagram will be extensively employed throughout this article, explicitly or implicitly. 

Observe also that since $\DD_{L'}\neq 0$ only if $L'$ consists of $n$ distinct cells in the positive quadrant, we can forget all the terms in the sum \ref{pkop} but those relative to such diagrams.
\end{rem}

\begin{prop}\label{proe}
Let $L$ be a lattice diagram. Then for any integer $k\ge 1$ we have
\begin{equation}\label{ekop}
e_k(\partial X)\Delta_L(X,Y)=\sum_{1\le i_1<i_2<\cdots<i_k\le n}\epsilon(L,e_k(i_1,\dots,i_k;L))\Delta_{e_k(i_1,\dots,i_k;L)}(X,Y)
\end{equation}
where $e_k(i_1,\dots,i_k;L)$ is the lattice diagram obtained by replacing the biexponents $(p_{i_1},q_{i_1}),\dots,(p_{i_k},q_{i_k})$ by $(p_{i_1}-1,q_{i_1}),\dots,(p_{i_k}-1,q_{i_k})$ and where the coefficient $\epsilon(L,e_k(i_1,\dots,i_k;L))$ is a positive integer. 
\end{prop}

\noindent{\it Proof.}
The proof is almost the same as for the previous proposition. We write
\begin{equation}
e_k(X)=\sum_{1\le j_1<\cdots< j_k\le n}x_{j_1}\dots x_{j_k}.
\end{equation}
We develop the determinantal form of $\Delta_L$ with respect to the columns $j_1,\dots,j_k$ to obtain the following expression where $\Delta_L^{i_1,\dots,i_k}$ denotes the lattice diagram determinant relative to the biexponents $i_1,\dots,i_k$ of $L$ and $A_{i_1,\dots,i_k;j_1,\dots,j_k}$ the cofactor:
\begin{equation}\label{veto}
\Delta_L=\sum_{1\le i_1<\cdots< i_k\le n}\Delta_L^{i_1,\dots,i_k}(x_{j_1},\dots ,x_{j_k}) A_{i_1,\dots,i_k;j_1,\dots,j_k}.
\end{equation}
Next we derive \ref{veto} to obtain 
\begin{eqnarray}
\partial (x_{j_1}\dots x_{j_k})\Delta_L=\sum_{1\le i_1<\cdots< i_k\le n}&\big(c_{i_1,\dots,i_k;j_1,\dots,j_k}\Delta_{e_k(i_1,\dots,i_k;L)}^{i_1,\dots,i_k}(x_{j_1},\dots ,x_{j_k})\nonumber \\
&\times A_{i_1,\dots,i_k;j_1,\dots,j_k}\big),
\end{eqnarray}
where $c_{i_1,\dots,i_k;j_1,\dots,j_k}$ is a positive integer. We see that $c_{i_1,\dots,i_k;j_1,\dots,j_k}$ is equal to $p_{i_1}\cdots p_{i_k}$ and thus does not depend on $j_1,\dots,j_k$ . Therefore we can omit the subscript $j_1,\dots,j_k$.
Thus we get
\begin{eqnarray}\label{reco}
e_k(\partial X)\Delta_L=\sum_{(1\le i_1< \cdots < i_k\le n)}\sum_{(1\le j_1<\cdots< j_k\le n)}&\!\!\!\!\!\!\!\! \big(c_{i_1,\dots,i_k}\Delta_{e_k(i_1,\dots,i_k;L)}^{i_1,\dots,i_k}(x_{j_1},\dots ,x_{j_k})\nonumber\\
&\times A_{i_1,\dots,i_k;j_1,\dots,j_k}\big).
\end{eqnarray}
By recognizing in \ref{reco} the development of $\Delta_{e_k(i_1,\dots,i_k;L)}$, we finally obtain the expected formula. The sign in front of the coefficient $\epsilon(e_k;i_1,\dots,i_k;L)$ should be the sign of the permutation that reorders the obtained biexponents in increasing lexicographic order. In fact the choice of the lexicographic order \ref{lex} is such that this permutation is always the identity: each cell stays in its original column and no one of them ``jumps'' over another one, so that the order is unchanged.

\begin{rem}\label{bonrem}\rm
A useful observation is the fact that in the Propositions \ref{proP} and \ref{proe}, the coefficient $\epsilon(L,L')$ only depends on the original diagram $L$ and on the final diagram $L'$, but not on the differential operator. Let us clearly define this coefficient: if $L=\{(p_1,q_1),\dots,(p_n,q_n)\}$ and $L'=\{(p'_1,q'_1),\dots,(p'_n,q'_n)\}$, $\epsilon$ is given by the following formula:
\begin{equation}\label{epsi}
\epsilon(L,L')=\frac{\prod_{i=1}^np_i!q_i!}{\prod_{i=1}^np'_i!q'_i!}.
\end{equation}
This coefficient is a positive integer that appears (up to sign) as the coefficient of $\DD_{L'}$ in the expression of $P(\partial X)\DD_L$, where $P$ is a power sum or an elementary symmetric function ; we shall see in the next proposition that it is also the case for homogeneous symmetric functions.

Another important remark is that we have to be careful when we apply products of differential operators. Indeed in this case multiplicities may appear in the formulas. Let $P(\partial)$ and $Q(\partial)$ be two differential operators such that formulas like \ref{pkop} or \ref{ekop} hold for $P(\partial)$ and $Q(\partial)$ with $\epsilon$ given by \ref{epsi}. We first observe that $\epsilon$ is multiplicative, i.e.
\begin{equation}
\epsilon(L,L')=\epsilon(L,L'')\epsilon(L'',L'),
\end{equation}
for $L$, $L''$ and $L'$ three diagrams.
Thus the coefficient of $\DD_{L'}$ in $P(\partial)Q(\partial)\DD_L$ is a multiple (up to sign for power sums) of $\epsilon(L,L')$.
This multiplicity corresponds to the number of choices in the order of the different shifts, that is to say the number of diagrams $L''$ such that $L''$ appears in $Q(\partial)\DD_L$ and $L'$ appears in $P(\partial)\DD_{L''}$.
This multiplicity is denoted by $c_{P,Q}(L'L')$

Let us take an example: if we apply $e_1(\partial X)e_1(\partial X)$ to the determinant of the diagram $L=\{(1,0),(1,1)\}$, we obtain a single diagram $L'=\{(0,0),(0,1)\}$, with $\epsilon(L,L')=1$, but
\begin{equation}
e_1(\partial X)e_1(\partial X)\DD_L=2\DD_{L'}.
\end{equation}
The multiplicity 2 correponds to the fact that we can either first move down the cell $(1,0)$ and next the cell $(1,1)$ or do it in the reverse order.

All these observations are crucial to well understand the proof of the following proposition.
\end{rem}

\vskip 0.2 cm
Now, to state the next proposition, we need to introduce some notation. For a lattice diagram $L$, we denote by $\overline L$ its complement in the
positive quadrant
(it is an infinite subset). Again we order $\overline L=\{(\overline
p_1,\overline q_1), (\overline
p_2,\overline q_2),\dots\,\}$ using the lexicographic order \ref{lex}.

\begin{prop}\label{h}
Let $L$ be a lattice diagram. Then for any integer $k\ge 1$ we have
\begin{equation}\label{eqh}
h_k(\partial X)\Delta_L(X,Y)=\sum_{1\le i_1<i_2<\cdots<i_k}
  \epsilon(L,h_k(i_1,\dots,i_k;L))\Delta_{h_k(i_1,\dots,i_k;L)}(X,Y)
\end{equation}
where
$h_k(i_1,\dots,i_k;L)$ is the lattice diagram with the following
complement diagram. Replace the
biexponents
$(\overline p_{i_1},\overline q_{i_1}),\dots,$ $(\overline
p_{i_k},\overline q_{i_k})$ of the complement
$\overline L$ with
$(\overline p_{i_1}+1,\overline q_{i_1}),\dots,(\overline
p_{i_k}+1,\overline q_{i_k})$ and keep the other unchanged. The coefficient
$\epsilon(L,h_k(i_1,\dots,i_k;L))$ is a positive integer, given by formula \ref{epsi}.
\end{prop}

\noindent{\it Proof.}
We shall prove this proposition by induction on $k$. If $k=1$, then $h_1=e_1$ and the result is true since moving down a cell is equivalent to moving up a hole. Assume the result is true up to $k-1$. Then we use the fact that $h_k=e_1h_{k-1}-e_2h_{k-2}+\cdots+(-1)^ke_{k-1}h_1+(-1)^{k+1}e_k$. 

Each term $e_{l}h_{k-l}$ for $1\le l\le k$ gives a linear combination of $\DD_{L'}$, whose coefficients are multiple of $\epsilon(L,L')$ according to Remark \ref{bonrem}. The problem is to compute the alternating sum of all these coefficients to get the result of $h_k(\partial X)\Delta_L$. 

Let $L'$ be one of the diagrams created by the terms $e_{l}h_{k-l}$. The coefficient of $\DD_{L'}$ in $e_{l}(\partial X)h_{k-l}(\partial X)$ is equal to $c_{e_l,h_{k-l}}(L,L')\epsilon(L,L')$. In this proof let us denote $c_{e_l,h_{k-l}}(L,L')$ simply by $c_{l}(L,L')$. The question is to compute 
\begin{equation}\label{somc}
\sum_{1\le l\le k}(-1)^{l+1}c_{l}(L,L'). 
\end{equation}
Let $k'\le k$ be the number of distinct holes moving between $L$ and $L'$ and $d\le k'$ the number of those which have a moving hole below them. Each of these $d$ holes has to move with $h_{k-l}(\partial X)$ because the hole below it is able to move up with $e_{l}(\partial X)$ only if it has a cell above itself. The choice therefore comes from the $k'-d$ other holes which can either move up with $h_{k-l}(\partial X)$ or not: we choose $k-l-d$ among them to move with $h_{k-l}(\partial X)$. Thus we get
\begin{equation}
c_{l}(L,L')={{k'-d}\choose{k-l-d}}={{k'-d}\choose{l-(k-k')}}.
\end{equation}
And the sum in \ref{somc} becomes
\begin{equation}
\sum_{l=1}^k(-1)^{l+1}{{k'-d}\choose{l-(k-k')}}=
\left\{
\begin{array}{cc}
1 & {\rm if}\ k'=k,\\
0 & {\rm if}\ k'<k.
\end{array}
\right. 
\end{equation}
Thus we get the desired formula \ref{eqh}.

\begin{rem}\rm
One efficient application of the previous proposition is to give a necessary condition that tests if a partial symmetric operator belongs to the vanishing ideal of a lattice diagram determinant (see \cite{nantelmoi}).
An example of the strength of this principle is to give immediate proofs of Propositions 1-2-3-4 of \cite{aval} (these propositions provide a Groebner basis of the vanishing ideal of $\Delta_{\mu}$ when $\mu$ is a hook). The previous proofs in \cite{aval} were recursive and intricate but the results now become simple applications of Proposition \ref{h}.
\end{rem}

\begin{rem}\rm
The shift operators are also useful to reduce the sum \ref{sumdef} defining $M_{i,j}^k$. In the special case of one hole, it is now easy to see that 
\begin{equation}\label{onehole}
M_{i,j}^1=M_{\mu/i,j}. 
\end{equation}
Indeed we have that for any integer $k$ and $l$
\begin{equation}
e_k(\partial X)e_l(\partial Y)\DD_{\mu/i,j}=c.\DD_{\mu/i+k,j+l},
\end{equation}
with $c$ an integer different from zero. This implies $M_{i,j}^1\subseteq M_{\mu/i,j}$, and the reverse inclusion is obvious.

In the particular case of two holes, let $k$ and $l$ be positive integers and let us use the following notations: for two cells $h_1$ and $h_2$, $\epsilon^i_{h_1,h_2}=\epsilon\big(\mu/\{(i,j),(i+1,j)\},\mu/\{h_1,h_2\}\big)$ and $\epsilon^j_{h_1,h_2}=\epsilon\big(\mu/\{(i,j),(i,j+1)\},\mu/\{h_1,h_2\}\big)$. If we are careful of the different signs by applying Propositions \ref{proP} and \ref{proe}, then we get the following identities
$$P_{l}(\partial Y)e_{k-1}(\partial X)\Delta_{\mu/\{(i,j),(i+1,j)\}}=P_{l}(\partial Y)\big(\epsilon^i_{(i,j),(i+k,j)} \DD_{\mu/\{(i,j),(i+k,j)\}}\big)$$
\vskip -0.4 cm
$$=(-1)^{b+h}\epsilon^i_{(i+k,j),(i,j+l)}\Delta_{\mu/\{(i+k,j),(i,j+l)\}}\ \ \ \ \ \ \ \ \ \ \ \ \ \ \ \ \ \ \ \ \ \ \ \ $$
\vskip -0.6 cm
\begin{equation}\label{tohol1}
\ \ \ \ \ \ \ \ \ \ \ \ \ \ \ \ \ \ \ \ \ \ \ \ +(-1)^{b+v+1}\epsilon^i_{(i,j),(i+k,j+l)}\Delta_{\mu/\{(i,j),(i+k,j+l)\}}
\end{equation}
\vskip 0.2 cm
\noindent
and
$$P_{k}(\partial\! X)e_{l-1}(\partial Y)\Delta_{\mu/\{(i,j),(i,j+1)\}}\!\!=\!\!P_{k}(\partial\! X)\big((-\!1)^{b+h+1}\epsilon^j_{(i,j),(i,j+l)} \DD_{\mu/\{(i,j),(i,j+l)\}}\big)$$
\vskip -0.6 cm
$$=(-1)^{b+h+1}\big((-1)^h\epsilon^j_{(i+k,j),(i,j+l)} \Delta_{\mu/\{(i+k,j),(i,j+l)\}}\ \ \ \ \ \ \ \ \ \ \ \ \ \ \ \ \ \ \ \ \ \ \ \ $$
\vskip -0.6 cm
\begin{equation}\label{tohol2}
\ \ \ \ \ \ \ \ \ \ \ \ \ \ \ \ \ \ \ \ \ \ \ \ +(-1)^v\epsilon^j_{(i,j),(i+k,j+l)}\Delta_{\mu/\{(i,j),(i+k,j+l)\}}\big),
\end{equation}

\begin{figure}[ht]
$$
\epsfbox{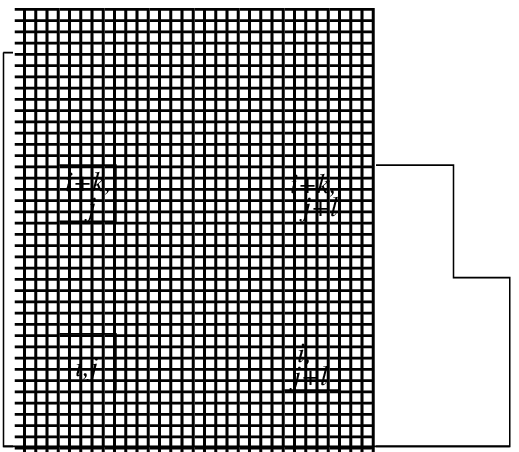}
$$
\end{figure}

\noindent
where $h$, $v$ and $b$ are respectively the numbers of cells with horizontal, vertical and both horizontal and vertical stripes in the figure above (we have to compute the sign of the permutation which reorders the cells in the lexicographic order \ref{lex}).

By observing that the product of the signs of the four coefficients in \ref{tohol1} and \ref{tohol2} is $(-1)^{2(2b+2h+v+1)+1}=(-1)$ we have that exactly three coefficients in \ref{tohol1} and \ref{tohol2} are of the same sign, whence $\Delta_{\mu/\{(i,j),(i+k,j+l)\}}$ and $\Delta_{\mu/\{(i+k,j),(i,j+l)\}}$ are in $M_{\mu/\{(i,j),(i,j+1)\}}+M_{\mu/\{(i,j),(i+1,j)\}}$.

Next, by Proposition \ref{h} we can move simultaneously the two holes.
This implies that for any couple of holes $(h_1,h_2)$ in the shadow of $(i,j)$ then $\Delta_{\mu/\{h_1,h_2\}}\in M_{\mu/\{(i,j),(i,j+1)\}}+M_{\mu/\{(i,j),(i+1,j)\}}$ thus 
\begin{equation}\label{twool}
M_{i,j}^2=M_{\mu/\{(i,j),(i,j+1)\}}+M_{\mu/\{(i,j),(i+1,j)\}}.
\end{equation}

The question of whether the obvious generalization of the previous result \ref{twool} is true when $k\ge 3$ appears naturally. Is it sufficient to take only the diagrams such that the holes form a partition of origin $(i,j)$? The answer is negative. For example it is easy to check (by computer) that when $\mu=(3,2)$,
\begin{equation}
\Delta_{\mu/\{(0,0),(1,0),(0,2)\}}\not \in M_{\mu/\{(0,0),(1,0),(0,1)\}}+M_{\mu/\{(0,0),(0,1),(0,2)\}}.
\end{equation}
\end{rem}

\section{The upper bound}

\begin{definition}
Let $M$ be a vector subspace of $\Q[Z]$ where $Z$ is a subalphabet of $(X_{n+k},Y_{n+k})$. We define its {\sl vanishing ideal} as the following ideal: 
\begin{equation}\label{kill}
I_M=\{P\in \Q[Z]\ :\ \forall Q\in M,\ P(\partial) Q=0\}. 
\end{equation}
If $P(\partial) Q=0$, we shall say that $P$ ``kills'' $Q$.

If $M={\mathcal L}_\partial[P]$ then we denote its vanishing ideal simply by $I_P$.
In the case of $M^k_{i,j}$ we denote $I_{M^k_{i,j}}$ by $I^k_{i,j}$.
\end{definition}
We recall the following important result (\cite{gh}, Proposition 1.1): 

\begin{prop}\label{perpe}
For $M$ a subspace of $\Q[Z]$, we have
\begin{equation}
M=I_M^{\perp}=\{P\in\Q[X_n,Y_n]:\ \forall Q\in I_M, \langle P,Q \rangle=0\}, 
\end{equation}
where the scalar product is defined by $\langle P,Q \rangle=L_0(P(\partial) Q)$ and where $L_0$ is the linear form that associates to a polynomial its term of degree 0.
\end{prop}

\subsection{About ideals}

We want here to prove the following

\begin{prop}\label{I=I}
\begin{equation}
I^k_{i,j}=\bigcap_{(a_1,b_1),\dots,(a_k,b_k)}I_{\partial x_{n+1}^{a_1}\partial y_{n+1}^{b_1}\cdots\partial x_{n+k}^{a_k}\partial y_{n+k}^{b_k}\Delta_{\mu}} \cap \Q[X_n,Y_n] \stackrel {\rm def} {=} {\I},
\end{equation}
where the intersection is over the $k$-tuples of different cells in the shadow of $(i,j)$ that we assume to be ordered in lexicographic order. 
\end{prop}

\noindent{\it Proof.}
Let $(a_1,b_1),\dots(a_k,b_k)$ be $k$ cells in $S_{\mu}((i,j))$, the shadow of $(i,j)$ in $\mu$.
By expanding $\Delta_{\mu}$ with respect to the last $k$ columns, we obtain:
$$\Delta_{\mu}(X_{n+k},Y_{n+k})=\sum_{(a'_1,b'_1),\dots,(a'_k,b'_k)}\pm\Delta_{\{(a'_1,b'_1),\dots,(a'_k,b'_k)\}}(\bar X_n,\bar Y_n)\ \ \ \ \ \ \ \ \ \ \ $$
\vskip -1 cm
\begin{equation}
\ \ \ \ \ \ \ \ \ \ \ \ \ \ \ \ \ \ \ \ \ \ \ \ \ \ \ \ \ \ \ \ \ \ \ \ \ \ \ \ \ \ \ \ \ \times\Delta_{\mu/\{(a'_1,b'_1),\dots,(a'_k,b'_k)\}}(X_n,Y_n),
\end{equation}
where $\bar X_n=\{x_{n+1},\dots,x_{n+k}\}$ and $\bar Y_n=\{y_{n+1},\dots,y_{n+k}\}$.
Thus we get:
\begin{equation}\label{etoile}
\partial (x_{n+1}^{a_1}y_{n+1}^{b_1}\cdots x_{n+k}^{a_k}y_{n+k}^{b_k})\Delta_{\mu}(X_{n+k},Y_{n+k})\! = \!c\Delta_{\mu/\{(a_1,b_1),\dots,(a_k,b_k)\}}(X_n,Y_n)+C
\end{equation}
where $c$ is a rational constant (different from 0) and $C$ a linear combination with coefficients in $\Q[x_{n+1},y_{n+1},\dots,x_{n+k},y_{n+k}]$ of polynomials $\Delta_{\mu/\{(a'_1,b'_1),\dots,(a'_k,b'_k)\}}\\(X_n,Y_n)$, with: 
\begin{equation}\label{sha}
\forall\ 1 \leq l \leq k,\quad (a_l^{\prime}, b_l^{\prime})\in S_{\mu}((i, j)). 
\end{equation}

Indeed $\Delta_{\{(a_1^{\prime}, b_1^{\prime}),
\dots,(a_k^{\prime}, b_k^{\prime})\}}(\bar X_n, \bar Y_n)$
is not killed by $\partial(x_{n + 1}^{a_1} y_{n + 1}^{b_1}
\cdots x_{n + k}^{a_k} y_{n + k}^{b_k})$ only if there
exists at least a permutation $\sigma \in {\cal S}_k$, the
symmetric group on $k$ elements, such that 
\begin{equation}\label{star}
(a_{\sigma(l)}^{\prime}, b_{\sigma(l)}^{\prime}) \in
S_{\mu}((a_l, b_l)), \quad \forall\ 1 \leq l \leq
k.
\end{equation} 
This follows easily from the definition of the
$\Delta_{\{(a_1^{\prime}, b_1^{\prime}), (a_2^{\prime},
b_2^{\prime}), \dots, (a_k^{\prime}, b_k^{\prime})\}}$
as a determinant: 
$$\Delta_{\{(a_1^{\prime}, b_1^{\prime}), (a_2^{\prime},
b_2^{\prime}), \dots, (a_k^{\prime}, b_k^{\prime})\}} =
\sum_{\sigma \in {\cal S}_k} {\rm sgn}(\sigma)\ x_{n +
1}^{a_{\sigma(1)}^{\prime}} y_{n
+ 1}^{b_{\sigma(1)}^{\prime}} x_{n +
2}^{a_{\sigma(2)}^{\prime}} y_{n
+ 2}^{b_{\sigma(2)}^{\prime}} \cdots x_{n +
k}^{a_{\sigma(k)}^{\prime}} y_{n
+ k}^{b_{\sigma(k)}^{\prime}}.$$ 
Taking the partial derivative $\partial(x_{n + 1}^{a_1} y_{n
+ 1}^{b_1} \cdots x_{n + k}^{a_k} y_{n + k}^{b_k})$, we
get \ref{star}. For all $1 \leq l \leq k$, we have  $S_{\mu}(a_l,
b_l) \subseteq S_{\mu}(i, j)$. Consequently 
$(a_{\sigma(l)}^{\prime}, b_{\sigma(l)}^{\prime}) \in
S_{\mu}(i, j), \ \forall\ 1 \leq l \leq k$. Because
$\sigma$ is a permutation, \ref{sha} is now obvious.    

To illustrate the equation \ref{etoile}, we give the following example: $\mu=(3,2)$, $n=3$, $k=2$, $(a_1,b_1)=(0,0)$, $(a_2,b_2)=(1,0)$, then
\begin{align*}
\partial(x_4^0y_4^0x_5^1y_5^0)\DD_{\mu}(X_5,Y_5)&=\DD_{\mu/\{(0,0)(1,0)\}}(X_3,Y_3)+y_5\DD_{\mu/\{(0,0)(1,1)\}}(X_3,Y_3)\\
&-y_4\DD_{\mu/\{(1,0)(0,1)\}}(X_3,Y_3)\\
&+(-x_4y_5+x_4y_4)\DD_{\mu/\{(1,0)(1,1)\}}(X_3,Y_3)\\
&-y_4^2\DD_{\mu/\{(1,0)(0,2)\}}(X_3,Y_3)+y_4y_5\DD_{\mu/\{(0,1)(1,1)\}}(X_3,Y_3)\\
&-y_4^2y_5\DD_{\mu/\{(1,1)(0,2)\}}(X_3,Y_3).
\end{align*}

\vskip 0.2 cm
\noindent
Hence we get what we want because:
\begin{itemize}
\item $\I\subseteq I_{i,j}^k$: let $P$ be a polynomial in $\I$. Since $P$ kills the left-hand side of \ref{etoile}, it kills the constant term in $\Q[\bar X_n,\bar Y_n]$ of the left-hand side which is $\Delta_{\mu/\{(a_1,b_1),\dots,(a_k,b_k)\}}(X_n,Y_n)$. Thus $P$ is in $I_{i,j}^k$.
\item $I_{i,j}^k\subseteq \I$: let $P$ be a polynomial in $I_{i,j}^k$. By \ref{sha}, $P$ kills all the terms of the right-hand side of \ref{etoile}; thus it kills the left-hand side. This implies $P\in\I$.
\end{itemize}

\subsection{Sets of points and vanishing ideals}

The reasoning is inspired from \cite{lattice}, Theorem 4.2.

Let $\mu$ be a partition of $n+k$, $l=\mu_1$ its length and $h$ its height (the number of its positive parts). We consider two sets $\alpha=(\alpha_0,\dots,\alpha_{h-1})$ and $\beta=(\beta_0,\dots,\beta_{l-1})$ of distinct rational numbers. To any injective tableau $T$ of shape $\mu$ with entries $\{1,\dots,n+k\}$, we associate a point $(a(T),b(T))$ in $\Q^{2(n+k)}$ by the following process: 
\begin{equation}\label{ab}
\forall\ 1\le t\le n+k,\ \ a_t(T)=\alpha_{r_t(T)}\ {\rm and}\ b_t(T)=\beta_{c_t(T)} 
\end{equation}
where $r_t(T)$ (resp. $c_t(T)$) is the number of the row (resp. column) of $T$ where the entry $t$ lies in $T$. We think useful to recall here that the convention introduced in the Definition \ref{mudef} is that the first row and column are indexed by 0. We define $\rho$ as the orbit of $(a,b)$ when $T$ varies over the $(n+k)!$ injective tableaux of shape $\mu$. Let us observe that, since the $\alpha_t$'s and $\beta_t$'s are distinct, two different tableaux give two different points, i.e. $T\mapsto(a(T),b(T))$ is an injective map.
We introduce $J_{\rho}$ the ideal of polynomials that are zero over all the orbit.
We recall that the operator $gr$ is the operator that associates to a polynomial its term of maximum degree and that the $gr$ of an ideal is the ideal generated by the $gr$ of its elements. Next we define $I=grJ_{\rho}$ and $H=I^{\perp}$.
We will use the following important result (cf. \cite{gh}, Theorem 1.1): 
\begin{prop}\label{ghmu}
For any choice of $\alpha_t$'s and $\beta_t$'s, if $I$ is the graded ideal associated to the vanishing ideal of $\rho$ then we have the inclusion
\begin{equation}
I\subseteq I_{\Delta_{\mu}}.
\end{equation}
\end{prop}
We now look at another set, this time in $\Q^{2n}$. We consider the set of tableaux $T$ of shape $\mu$ with $n$ entries and $k$ white cells such that the $k$ white cells are in the shadow of $(i,j)$. Let us denote this set of tableaux by ${\T}_{i,j}^k$.
By the same process as described in \ref{ab}, we define a set $\rho^k$ in $\Q^{2n}$.
We recall here that we denote by $s_{\mu}(i,j)$ or simply by $s$ the cardinality of the shadow of the cell $(i,j)$ in $\mu$.
Since the cardinality of ${\T}_{i,j}^k$ is ${s \choose k} n!$ and the process \ref{ab} is still injective, the set $\rho^k$ has ${s \choose k} n!$ points. 
We introduce $J_{\rho^k}$ the ideal of polynomials that are zero all over $\rho^k$, and $I^k=grJ_{\rho^k}$ and $H^k=(I^k)^{\perp}$. 

The first information is given by the following equation: 
\begin{equation}\label{dihk}
\dim H^k={s \choose k} n!.
\end{equation}
This comes from the fact that $\dim H^k=\dim \Q[X,Y]/J_{\rho^k}=\#\rho^k={s \choose k} n!$. We shall not develop this point, extensively treated in \cite{orbit}.  

We want to prove that $M_{i,j}^k\subseteq H^k$ and by Proposition \ref{perpe}, it is equivalent to prove that $I^k\subseteq {\I}.$

\subsection{Inclusion}

We want here to obtain the next proposition:

\begin{prop}

We have the inclusion:
\begin{equation}
I^k\subseteq {\I}.
\end{equation}
\end{prop}

\noindent{\it Proof.}
Let $P$ be a polynomial in $J_{\rho^k}$. Let us consider
\begin{eqnarray}
Q(X_{n+k},Y_{n+k})=&\!\!\!\!\!P(X_n,Y_n)\times\prod_{i'=0}^{i-1}(x_{n+1}-\alpha_{i'})\cdots\prod_{i'=0}^{i-1}(x_{n+k}-\alpha_{i'})\nonumber\\
&\times\prod_{j'=0}^{j-1}(y_{n+1}-\alpha_{j'})\cdots\prod_{j'=0}^{j-1}(y_{n+k}-\alpha_{j'}).
\end{eqnarray}
We want to check that this polynomial is an element of $J_{\rho}$. We take an element $(a,b)=(a(T),b(T))$ of $\rho$. If its projection on $\Q^{2n}$ (by keeping the first $n$ entries of $a$ and $b$) is in $\rho^k$ then $Q(a,b)=0$ because of $P$. If not, the tableau $T$ must have at least one entry between $n+1$ and $n+k$ in the complement of the shadow of $(i,j)$, i.e. in the first $i$ rows or the first $j$ columns and we have still $Q(\alpha,\beta)=0$.

Thus $Q\in J_{\rho}$, hence $gr(Q)\in I_{\Delta_{\mu}}$. Next by looking at the term of maximal degree we get: $gr(P)\in I_{\partial x_{n+1}^{i+1}\partial y_{n+1}^{j+1}\cdots\partial x_{n+k}^{i+1}\partial y_{n+k}^{j+1}\Delta_{\mu}}$.

For any set of k cells $\{(a_1,b_1),\dots,(a_k,b_k)\}$ in the shadow of $(i,j)$, we observe that $\forall r,\  1\le r\le k$, $a_r\ge i$ and $b_r\ge j$. 
Hence $gr(P)$ is in $\I$, which was to be proved.

\subsection{Conclusion}

The main result is now a consequence of all what precedes:

\begin{thm}\label{ub}
If $\mu$ is a partition of $n+k$ and $s$ the cardinality of the shadow of the cell $(i,j)$, then we have:
\begin{equation}
\dim M^k_{i,j}\le {s \choose k} n!.
\end{equation}
\end{thm}

\begin{rem}\rm
If we recall the proof of Theorem 1.1 of \cite{gh}, we observe that the previous reasoning implies the following fact. If equality holds in Theorem \ref{ub}, then $M^k_{i,j}$ decomposes as ${s \choose k}$ times the left regular representation.
\end{rem}

\vskip 0.2 cm

Numerical examples and the fact that the construction described in the previous subsection affords the ``good'' upper bound in the case of one set of variables (see the next section) support the following conjecture, which was first stated by F. Bergeron.

\begin{conj}\label{conje}
{\it
With the notations of the previous theorem:
\begin{equation}
\dim M^k_{i,j} = {s \choose k} n!.
\end{equation}
}
\end{conj}

\begin{rem}\rm
When $k=1$, this conjecture reduces to Conjecture I.2 of \cite{lattice} and when $s=k$ or $k=0$ to the $n!$ conjecture.
\end{rem}

\section{Case of one set of variables}

\begin{definition}
Let $M=M(X,Y)$ be a subspace of $\Q[X,Y]$. Then we denote by $M(X)$ the subspace of $M$ consisting of elements of 0 $Y$-degree. We also denote the vanishing ideal of $M(X)$ by $I_M(X)$.
\end{definition}

The goal of this section is to obtain an explicit basis for $M_{i,j}^k(X)$, the subspace of $M_{i,j}^k(X,Y)$ of elements of 0 $Y$-degree.

\subsection{Construction}

We first recall results about $M_{\mu}(X)$ the subspace of $M_{\mu}$ of elements of 0 $Y$-degree (which was denoted by $M_{\mu}^0$ in \cite{aval} and \cite{untrou}). When $\mu$ is a partition of $n$, we have
\begin{equation}\label{muze}
\dim M_{\mu}(X)=n!/\mu!, 
\end{equation}
where $\mu!=\mu_1!\cdots\mu_k!$. This space has been studied in \cite{aval}, \cite{gh}, \cite{nantel adriano}. 
Let $\M(\mu)$ be a set of monomials whose cardinality is $n!/\mu!$, such that the set $B_{\mu}=\{M(\partial)\Delta_{\mu}:\ M\in \M(\mu)\}$ is a basis for the space $M_{\mu}(X)$.
By the work in \cite{aval} we know such a set exists.

Now let $\mu$ be a partition of $n+k$. Next we choose in the Ferrers diagram $\mu$, $k$ cells which are simultaneously in the shadow of $(i,j)$ and such that any circled cell has either a cell outside the partition on its right or a circled cell (see the figure below). A circled cell satisfying this condition is said to be ``Right''. We denote by $\F_{\mu}^k$ the set of the obtained objects, which we call Right diagrams (associated to $\mu$). 

\begin{figure}[ht]\label{cons}
\caption{A Right diagram $F$ and its associated partition $\mu_F$ and diagram with $k$ holes $\mu_F^k$.}
$$
\epsfbox{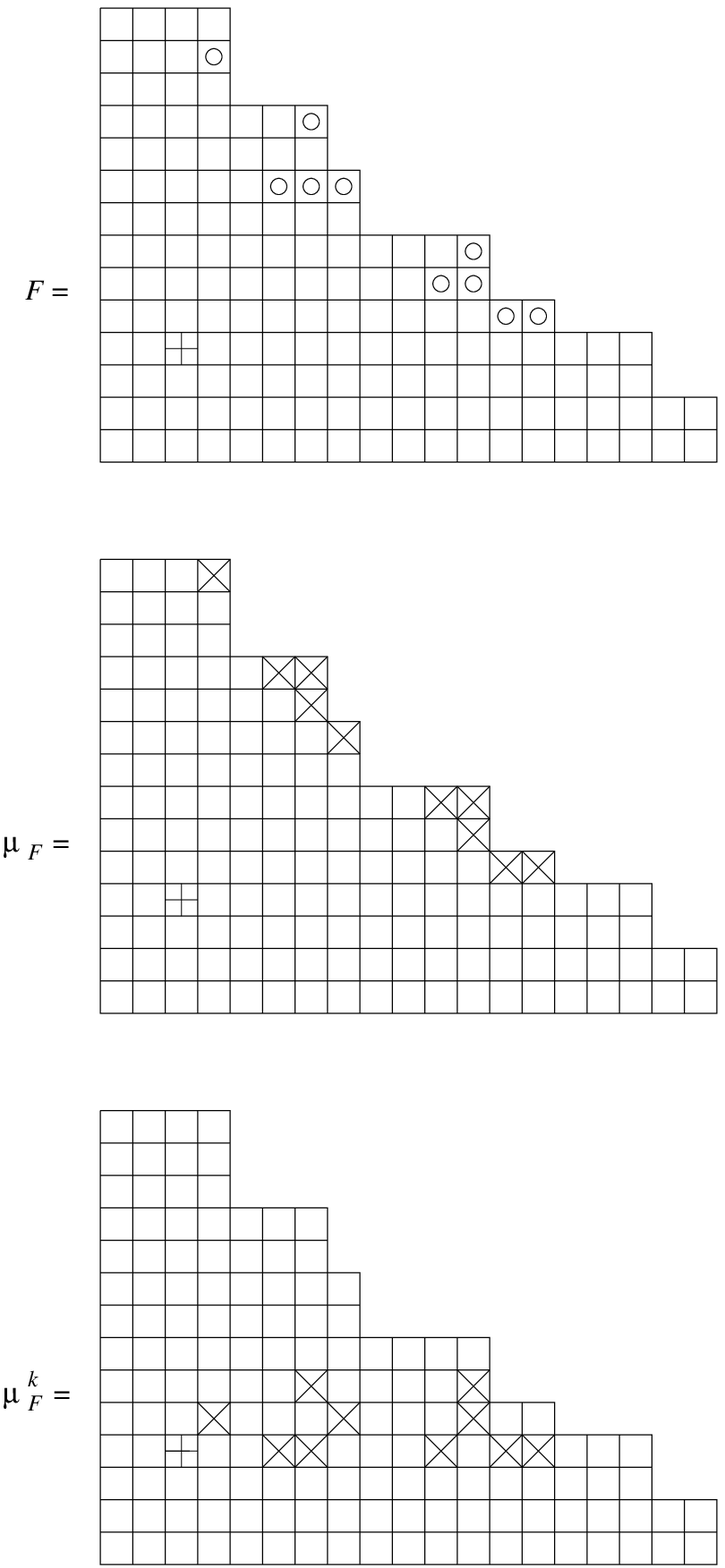}
$$
\end{figure}

We are now going to associate to each Right diagram two objects: a partition and a diagram with (at most) $k$ holes in the shadow of $(i,j)$.
The Figure 1 illustrates this construction.
In this figure, the chosen cells in the Right diagram are cells with a circle, in the cell $(i,j)$ appears a $+$ sign and the holes are as usual cells with crosses ($\times$). In this example $n=142$ and $k=10$.

To a Right diagram $F$ in $\F_{\mu}^k$ we first associate $\mu_{F}$ the partition of $n$ obtained by pushing up the circled cells and by removing the corresponding cells (see Figure 1). 

We also define a diagram $\mu_F^k$ with $k$ holes by proceeding as follows. We look at the columns where a circled cell appears. In our example we have $8$ such columns. For a column $j'\ge j$ where a circled cell appears, we denote by $h(j')$ the number of places where we could have put a circled cell (of course a Right one) below the lowest circled cell of this column. In our example, we have: $h(3)=1,\ h(5)=0,\ h(6)=0,\ h(7)=1,\dots,\ h(13)=0$. Next for any column $j'$ with a circled cell, we do the following. We denote the positions of the circled cells in this column by $(c(j'),j'),\ (c(j')+a_1,j'),\ \dots,\ (c(j')+a_d,j')$, with $(c(j'),j')$ the position of the lowest one, $0<a_1<\cdots<a_d$, and $d+1$ the number of circled cells in the column $j'$. We then place holes in cells $(i+h(j'),j'),\ (i+h(j')+a_1,j'),\ \dots,\ (i+h(j')+a_d,j')$. Doing this for all columns gives the diagram $\mu_F^k$.
This construction is illustrated in Figure 1. 

The crucial idea is to apply the monomials associated to $\mu_F$ to the determinant associated to $\mu_F^k$ and we are now able to state the main result of this section.

\begin{thm}\label{teo}
With the previous notations
\begin{equation}\label{base}
B_{i,j}^k(X)=\{M(\partial)\Delta_{\mu_F^k}:\ M\in \M(\mu_F),\ F\in \F_{\mu}^k\}
\end{equation}
is a basis for $M_{i,j}^k(X)$.
\end{thm}

The object of the end of the article is to prove this theorem. We will obtain an upper bound for the dimension of $M_{i,j}^k(X)$, next verify that the cardinality of $B_{i,j}^k(X)$ is equal to this upper bound, and prove that the family $B_{i,j}^k(X)$ is linearly independent.

\clearpage

\subsection{Upper bound}

\begin{definition}\label{deft}
We denote by $\T_{i,j}^k$ the set of injective, row-increasing tableaux of shape $\mu$ with $n$ entries $\{1,\dots,n\}$ and $k$ white cells (without any entry) such that the $k$ white cells are in the shadow of $(i,j)$ and not on the left side of an entry $1,\dots,n$.

We can also see $\T_{i,j}^k$ as the set of injective, row-increasing tableaux with entries $1,\dots,n$ of shapes all the Right diagrams $F$ of $\F_{\mu}^k$.
\end{definition}

The following lemma will be useful in the proof of the next proposition.
\begin{lem}\label{bonle}
Let $M=M(X,Y)$ be a subspace of $\Q[X,Y]$ and $M(X)$ its subspace of elements of 0 $Y$-degree. We suppose that $M$ is stable under derivation. Then we have the following relation between vanishing ideals:
\begin{equation}
I_M(X)=I_M\cap\Q[X].
\end{equation}
\end{lem}

\noindent{\it Proof.}
The inclusion $I\cap\Q[X]\subseteq I(X)$ is immediate. The reverse inclusion is obtained as follows. If $P$ is an element of $I(X)$ and $Q$ a polynomial in $M(X,Y)$, we look at the monomials of $Q$ in $Y$ with coefficients in $\Q[X]$. These coefficients are elements of $M(X)$ because $M$ is supposed to be stable under derivation. Thus these coefficients are killed by $P$ and so is $Q$ itself. 

\vskip 0.2 cm
The next proposition gives the analogue upper bound to Theorem \ref{ub} in the case of one set of variables.

\begin{prop}\label{upper bound}
The dimension of $M_{i,j}^k(X)$ satisfies the following inequality:
\begin{equation}\label{bsu}
\dim M_{i,j}^k(X)\le \#\T_{i,j}^k.
\end{equation}
\end{prop}

\noindent{\it Proof.}
From the Proposition \ref{I=I} and the Lemma \ref{bonle} applied to $M_{i,j}^k(X,Y)$, which is of course stable under derivation, we deduce that $I^k_{i,j}(X)=I^k_{i,j}\cap\Q[X_n]=\I\cap\Q[X_n]$.

We consider the projection of the set $\rho^k$ on $\Q^n$, i.e. we associate to each injective tableau $T$ of shape $\mu$ with $n$ entries $\{1,\dots,n\}$ and $k$ holes in the shadow of $(i,j)$ the point $a_{|n}(T)$ following the process defined in \ref{ab}. 
Let $\rho_0^k$ denote this set of points and $J_{\rho^k}^0$ its vanishing ideal. 
From the definition of $a_{|n}(T)$ it is clear that two tableaux give the same point if and only if they have the same entries on each line.
It is therefore equivalent to associate a point $a_{|n}(T)$ to each tableau $T$ in $\T_{i,j}^k$. In this case the correspondance is one-to-one and the number of points in $\rho_0^k$ is precisely $\#\T_{i,j}^k$, which is also the dimension of $gr(J_{\rho^k}^0)^{\perp}$ (this is the analogue of \ref{dihk}).

It remains to prove the following inclusion to justify Proposition \ref{upper bound}: 
\begin{equation}
gr(J_{\rho^k}^0) \subseteq I^k_{i,j}(X). 
\end{equation}
Let $P$ be a polynomial in $J_{\rho^k}^0$. Since $P\in\Q[X_n]\subset\Q[X_n,Y_n]$, $P$ is also in the vanishing ideal of $\rho^k$, thus $gr(P)\in I^k_{i,j}$ and next $gr(P)\in I^k_{i,j}\cap\Q[X_n]=I^k_{i,j}(X)$. Hence we have $gr(J_{\rho^k}^0) \subseteq I^k_{i,j}(X)$ and the equation \ref{bsu} is now a consequence of Proposition \ref{perpe}.

\subsection{Cardinality}

We claim that:
\begin{prop}
We have the following equality
\begin{equation}
\#B_{i,j}^k(X)=\#\T_{i,j}^k.
\end{equation}
\end{prop}

\noindent{\it Proof.}
Let $h$ be the height of the partition $\mu$. For a fixed Right diagram $F$ in $\F_{\mu}^k$, the number of associated elements in $B_{i,j}^k(X)$ is equal to $\frac {n!} {r_1!\cdots r_h!}$ where the $r_t$'s are the lengths of the rows of $\mu_F$ because of \ref{muze} and \ref{base}. By Definition \ref{deft} the number of elements in $\T_{i,j}^k$ associated to $F$ is $\frac {n!} {s_1!\cdots s_h!}$ where the $s_t$'s are the lengths of the rows of $F$.

It is therefore sufficient to observe that we do not change the cardinality by pushing up the circled cells. We look at the example of the lines 9, 10 and 11 of the previous example.
\vskip 0.3 cm

\centerline{
\epsffile{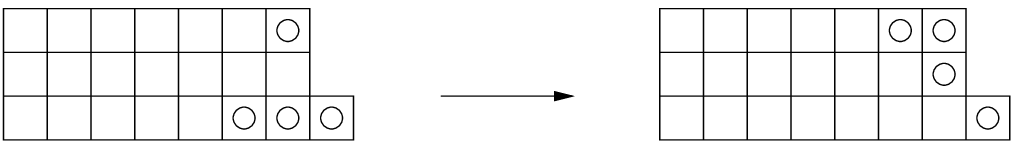}}

\vskip 0.3 cm
We observe that the lengths of the lines before the transformation are 5, 7, 6 and after the transformation 7, 6, 5. Thus the set of the lengths is unchanged. It is easy to see that it is always the case: the operation that pushes the holes up only permutes the lengths of the rows.

\subsection{Independence}

We want here to conclude the proof of the Theorem \ref{teo} by proving the independence of the set $B_{i,j}^k(X)$.

\begin{prop}\label{indep}
The set of polynomials $B_{i,j}^k(X)$ defined in Theorem \ref{teo} is linearly independent. Thus in particular equality holds in Proposition \ref{upper bound}.
\end{prop}

\noindent{\it Proof.}
Assume that we have a non-trivial dependence relation. 

We define the depth of a hole to be the number of cells (different from holes) that are above this hole.
We look at the $k$-tuples of the depths of the $k$ holes of $\mu_F^k$: $(d_1\le d_2\le\dots\le d_k)$. The crux of the proof is the following result:

\begin{lem}
The $k$-tuples $(d_1,d_2,\dots,d_k)$ are all distinct.
\end{lem}

\noindent{\it Proof.}
We want to prove that the depth of the holes increases from the right to the left and from top to bottom, and that two different Right diagrams $F$  and $F'$ of $\F_{\mu}^k$ give two different $k$-tuples of depths. We look at the circled cells with respect to this order.
We refer to the next figure and look at the columns from the right to the left. In this figure, $c$ denotes the number of circled cells in the column that we consider, $m$ the number of positions below the lowest circled cell where we could put a circle (these cells appear with a square), $l$ the height of this column (we look only at the cells above the $i$-th row) and $l+h$ is the height of the ``next'' column (i.e. the first on the left).
We want to prove that if we put a circled cell either in this column or in the next one, its depth will be greater or equal to the preceeding ones, and that its position is unambiguous if its depth is given.

\vskip 0.3 cm

\centerline{
\epsffile{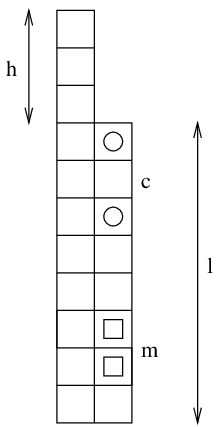}}

\vskip 0.3 cm

The depth of the lowest circled cell is $p=l-c-m$. The highest depth that could be obtained in this column is $l-c$ if $m=0$ and $l-c-1$ if $m>0$. In the next column the lowest depth is (it corresponds to put a circle at the top of the column): $l+h-1-c-h+1=l-c$. Thus there is no ambiguity for the position of the next circle if its depth is given, which proves the lemma.

\vskip 0.2 cm

Now let us complete the proof of Proposition \ref{indep}. If we have a non-trivial dependence relation between the elements of $B_{i,j}^k(X)$, we consider the greatest $k$-tuple of depths with respect to the lexicographic order which appears in this relation: $(d_1^0,d_2^0,\dots,d_k^0)$. This $k$-tuple is relative to a Right diagram $F^0$. We then apply the differential operator $h_k(\partial)^{d_1^0}.h_{k-1}(\partial)^{d_2^0-d_1^0}\dots h_1(\partial)^{d_k^0-d_{k-1}^0}$ to the dependence relation.
It kills all the terms but those which come from the single Right diagram $F^0$. These terms give polynomials which are in $B=\{M(\partial).\Delta_{\mu_{F^0}}:\ M\in \M(\mu_{F^0})\}$. They are independent since $B$ is a basis of $M_{\mu_{F^0}}(X)$.

The proof of Proposition \ref{indep} and as a consequence of Theorem \ref{teo} are now complete.

\vskip 0.5 cm
\noindent {\bf \large Acknowledgement.}
The author would like to thank François Bergeron for indicating him this problem and for numerous valuable suggestions and also the referees for their efforts to improve the quality of this article.

\end{document}